\documentclass[11pt]{article}%
\usepackage{amsfonts}
\usepackage{graphicx}
\usepackage{amsmath}
\usepackage{amssymb}
\usepackage{setspace}%
\providecommand{\U}[1]{\protect\rule{.1in}{.1in}}
\newtheorem{theorem}{Theorem}

\newtheorem{corollary}[theorem]{Corollary}

\newtheorem{lemma}[theorem]{Lemma}

\newtheorem{remark}[theorem]{Remark}

\singlespace
\begin{document}

\begin{center}
{\Large On kernel estimators of density for reversible Markov chains} \bigskip

\centerline{\today} \bigskip Martial Longla$^{a}$, Magda Peligrad$^{b}%
$\footnote{Supported in part by a Charles Phelps Taft Memorial Fund grant, and
the NSF grant DMS-1208237.} and Hailin Sang$^{a}$
\end{center}

\bigskip$^{a}$ Department of Mathematics, University of Mississippi,
University, MS 38677, USA. E-mail addresses: mlongla@olemiss.edu, sang@olemiss.edu

$^{b}$ Department of Mathematical Sciences, University of Cincinnati, PO Box
210025, Cincinnati, OH 45221-0025, USA. E-mail address: peligrm@ucmail.uc.edu

\bigskip



\begin{center}
\bigskip\textbf{Abstract}
\end{center}

In this paper we investigate the kernel estimator of the density for a
stationary reversible Markov chain. The proofs are based on a new central
limit theorem for a triangular array of reversible Markov chains obtained
under conditions imposed to covariances, which has interest in itself.

\section{Introduction and main results}

For estimating the marginal density for dependent sequences the dependence
structure plays an important role. One possible estimator is the kernel
estimator introduced by Rosenblatt (1956a). In general the dependence is
imposed in terms of mixing conditions (Bradley, 1993, Bosq et al, 1999 among
many others), in terms of coupling coefficients\ for functions of i.i.d. (Wu
et al, 2010) or positive association of random variables (Lin, 2003).

In this paper we study the kernel estimator for reversible Markov chains. It
is well known that for strictly stationary reversible Markov chains the
covariances can be viewed as a measure of dependence (see Kipnis and Varadhan,
1986). When estimating the density via kernel estimators we introduce a
triangular array of random variables which is only row-wise stationary. This
makes it difficult for studying the kernel density of the marginal
distribution for reversible Markov chains without imposing recurrence
conditions. As a matter of fact results on the kernel estimators for marginal
density of reversible Markov chains are very rare. We noticed only the paper
by Lei (2006) dealing with large deviations results for the integrated error
of the kernel density estimators for reversible Markov chains. The class they
considered is of reversible irreducible Markov chains with the transitions
satisfying a uniform integrability condition in square mean. However their
result cannot be applied when studying the density at a point or several
points. In this paper we develop tools that make this study possible.

Let $(X_{n})_{n\in\mathbb{Z}}$ be a stationary reversible Markov chain with
marginal distribution $\pi(A)=P(X_{n}\in A),$ for all Borel sets $A$. For a
stationary Markov chain the reversibility means that the distribution of
$(X_{0},X_{1})$ is the same as of $(X_{1},X_{0}).$ Assume that $\pi$ has a
marginal density $f(x)$, continuous at $x$. We shall consider in this paper
the Rosenblatt (1956a) estimator of density defined by
\begin{equation}
\hat{f}_{n}(x)=\frac{1}{nb_{n}}\sum_{k=1}^{n}K(\frac{x-X_{k}}{b_{n}}),
\label{kde}%
\end{equation}
where $b_{n}$ is a bandwidth converging to $0$ and $K$ is a kernel, a known
density function.

The problem considered in this paper is the consistency and the speed of
convergence of the kernel density estimator of the density at several points
$(x_{j})_{1\leq j\leq m}$ which will be given via a multivariate CLT.

To treat the problem we shall use the following notation
\begin{equation}
H_{k}(u,v)=P(X_{0}>u,X_{k}>v)-P(X_{0}>u)P(X_{k}>v) \label{defH}%
\end{equation}
and we shall denote
\begin{equation}
\eta_{k}=\int\int|H_{k}(u,v)|dudv. \label{def eta}%
\end{equation}
The condition we shall impose to $\eta_{k}$ is
\begin{equation}
\eta_{k}\leq\frac{1}{k^{4}l(k)}, \label{condeta}%
\end{equation}
where $l(x)\ $is a function increasing to infinity such that for any positive
$\ k,$ lim$_{x\rightarrow\infty}l(kx)/l(x)=1$ (slowly varying at infinite).

The following condition is imposed to the joint density of the vector
$(X_{0},X_{2})$ and a family of points of interest $(x_{j})_{1\leq j\leq m}$:
there exists the joint density $f_{2}(x,y)$ of $(X_{0},X_{2})$ which is
locally bounded around any pair $(x_{i},x_{j})_{1\leq i,j\leq m}$ in the sense
that there exists a constant $M$ and a constant $C_{M}$ (both depending on
$(x_{i},x_{j})$) such that
\begin{equation}
\sup_{|a|<M}|f_{2}(x_{i}+a,x_{j}+a)|<C_{M}. \label{second bound}%
\end{equation}
This condition is weaker than the condition which is usually imposed in the
dependent cases which requires that all the densities of vectors $(X_{0}%
,X_{j})$ are uniformly bounded on $\mathbb{R}^{2}$ (see condition in Bosq,
1998, or in Bosq et al, 1999). Local conditions can be found for instance in
papers by Liebscher (1999) and Dedecker and Merlev\`{e}de (2002).

All along the paper, we assume that the kernel $K$ satisfies the Condition C below:

(C1) $K$ is symmetric decreasing on $(0,\infty)$ and $\int K(u)du=1.$

(C2) $x^{2}K(x)\rightarrow0$ as $x\rightarrow\infty.$

(C3) $K$ is differentiable with $K^{\prime}(x)$ bounded.

Note that a normal kernel will satisfy all these conditions. The convergence
in distribution will be denoted by $\Rightarrow,$ and $\xrightarrow{P}$
denotes the convergence in probability.

The main result of this paper is the following:

\begin{theorem}
\label{gclt} Let $(X_{j})_{j\in\mathbb{Z}}$ be a stationary reversible Markov
chain with marginal density function $f(x)$ satisfying condition
(\ref{condeta}). Assume that the bandwidth $b_{n}$ in the estimator
(\ref{kde}) satisfies $nb_{n}^{4}\rightarrow\infty$ and the kernel $K$
satisfies Condition $C$. Then, at any points $x_{1},\cdots,x_{m}$ where $f(x)$
is continuous, different of $0$ and the joint densities satisfy condition
(\ref{second bound}), we have
\[
\sqrt{nb_{n}}\left(  \frac{\hat{f}_{n}(x_{j})-\mathbb{E}\hat{f}_{n}(x_{j}%
)}{(\hat{f}_{n}(x_{j})\int K^{2}(u)du)^{1/2}},1\leq j\leq m\right)
\Rightarrow N(0,I_{m}),
\]
where $I_{m}$ is the identity matrix.
\end{theorem}

It is well known that if the density is twice continuously differentiable at
$x_{j}$ then the bias is of order (see H\"{a}rdle 1991, relation (2.3.2))
\[
\mathbb{E}(\hat{f}_{n}(x_{j}))-f(x_{j})=\frac{b_{n}^{2}}{2}f^{\prime\prime
}(x_{j})+o(b_{n}^{2})\text{ as }b_{n}\rightarrow0.
\]
By combining this result with Theorem \ref{gclt} we get the following corollary:

\begin{corollary}
In addition to the conditions of Theorem \ref{gclt}, assume that $f$ is twice
continuously differentiable at $(x_{j})_{1\leq j\leq m}$ and $nb_{n}%
^{5}\rightarrow0.$ Then%
\[
\sqrt{nb_{n}}\left(  \frac{\hat{f}_{n}(x_{j})-f(x_{j})}{(\hat{f}_{n}%
(x_{j})\int K^{2}(u)du)^{1/2}},1\leq j\leq m\right)  \Rightarrow N(0,I_{m}).
\]

\end{corollary}

\bigskip

Let us comment about the dependence coefficient used in our results defined in
(\ref{def eta}).

If we have positive dependence, in the sense that $H_{k}(x,y)\geq0$ for all
$(x,y)\in\mathbb{R}^{2}$, then by the Lehmann Lemma (see Newman, 1980) we
have
\[
\eta_{k}=\mathrm{cov}(X_{0},X_{k}).
\]

This coefficient can also be controlled by a pairwise mixing condition which
is weaker than strong mixing coefficient introduced by Rosenblatt (1956b). As
in Rio (2000) relation (1.8a), define%
\[
\bar{\alpha}_{k}=\bar{\alpha}_{k}(X_{0},X_{k})=2\sup_{(x,y)\in\mathbb{R}^{2}%
}|H_{k}(x,y)|
\]
Theorem 1.1 in Rio (2000) states an estimate of the covariance between $X_{0}$
and $X_{k}$ in terms of $\alpha_{k}.$ However, in the proof, the author
actually estimated $\eta_{k}.$ Therefore we have
\[
\eta_{k}\leq2%
{\displaystyle\int\limits_{0}^{\bar{\alpha}_{k}}}
Q_{|X_{0}|}^{2}(u)du,
\]
where $Q_{|X_{0}|}$ is the quantile function of $|X_{0}|,$ i.e. the
generalized inverse of the function $\mathbb{P}(|X_{0}|>t).$

In particular, for $\delta>0$
\[
\eta_{k}\leq2\bar{\alpha}_{k}^{\delta/(2+\delta)}||X_{0}||_{2+\delta}^{2},
\]
and if $||X_{0}||_{\infty}\leq1$ then $\eta_{k}\leq2\bar{\alpha}_{k}$.

Recall that, in the stationary setting, the Rosenblatt pairwise strong mixing
coefficient, for an integer $k>0,$ is defined by%
\[
\alpha_{k}=\alpha_{k}(X_{0},X_{k})=\sup_{A,B\in\mathcal{B}}|\mathbb{P}%
(X_{0}\in A,X_{k}\in B)-\mathbb{P}(X_{0}\in A)(X_{k}\in B)|,
\]
where $\mathcal{B}$ denotes the Borel sigma algebra on the line. Clearly
$\bar{\alpha}_{k}\leq2\alpha_{k}.$ In terms of these mixing coefficients we
make the following remark.

\begin{remark}
\label{rmclt}Theorem \ref{gclt} also holds if we replace condition
(\ref{condeta}) by $\sum_{k\geq1}k\alpha_{k}<\infty.$
\end{remark}

It should be noted that this strong mixing rate was already pointed out in
Bosq et al. (1999), without assuming reversibility. The advantage here is that
we can have the asymptotic normality of the kernel estimators by only
requiring a condition on the joint density of the vector $(X_{0},X_{2})$ and
not on all the joint densities.

\bigskip

We finish this section by mentioning a few notations which will be used in
this paper. The largest integer smaller or equal to $x$ will be denoted by
$[x].$ By $c_{n}\ll d_{n}$ we understand that $c_{n}\leq Cd_{n}$ for some
$C>0$ and all $n;$ we denote by $a\vee b$ the maximum between $a$ and $b.$

\section{\textbf{Technical results}}

The proof relies on the properties of reversible Markov chains. We shall point
first some monotonicity conditions for integrals of functions of reversible
Markov chains. The regular conditional probability of $X_{1}$ given $X_{0}$
will be denoted by $Q(x,A)=\mathbb{P}(X_{1}\in A|\,X_{0}=x)$. Let $Q$ also
denote the Markov operator {acting via $(Qf)(x)=\int_{S}f(s)Q(x,ds).$ Next,
let $\mathbb{L}_{2}^{0}(\pi)$ be the set of measurable functions such that
$\int g^{2}d\pi<\infty$ and $\int gd\pi=0.$ In operator terms }the Markov
chain is called reversible if $Q=Q^{\ast},$ where $Q^{\ast}$ is the adjoint
operator of \ $Q$.

{For some function } ${g}\in$ {$\mathbb{L}_{2}^{0}(\pi)$, let }%
\[
{Y_{i}=g(X_{i}).}%
\]
Denote $\mathcal{F}_{n}=\sigma(\cdots X_{n-1},X_{n});$ $\ \mathbb{E}%
_{k}Y=\mathbb{E}(Y|\mathcal{F}_{k}).$ It is well known that, by the definition
of Markov chains, for $Y\in\sigma(X_{i},i\geq1)$ we have $\mathbb{E}%
_{0}Y=\mathbb{E}(Y|X_{0})=\mathbb{E}_{X_{0}}(Y).$ From the spectral theory of
self-adjoint operators on Hilbert spaces (see for instance Rudin, 1991), it is
well known that for every $g\in\mathbb{L}_{2}^{0}(\pi)$ there is a unique
\emph{transition spectral measure} $\nu$ supported on the spectrum of the
operator $[-1,1]$, such that
\begin{equation}
\mathbb{E}(\mathbb{E}_{0}(Y_{i})\mathbb{E}_{0}(Y_{j}))=\int_{-1}^{1}s^{i+j}%
\nu(\mathrm{d}s). \label{spectral}%
\end{equation}
By using this representation we give the following lemma which relates the
conditional expectation with the covariances for functions of reversible
Markov chains. It also points out several monotonicity conditions for the covariances.

\begin{lemma}
\label{cov}For stationary reversible Markov chains and every positive integers
$k,j$ we have
\begin{equation}
\mathbb{E}(\mathbb{E}_{0}Y_{k}\mathbb{E}_{0}Y_{j})=\mathbb{E}(Y_{0}Y_{k+j}).
\label{l1}%
\end{equation}
For any integer $k\geq0$%
\begin{equation}
\mathbb{E}(Y_{0}(Y_{2k}+Y_{2k+1}))\geq0, \label{l2}%
\end{equation}
also, for any integers $j$ and $k$ such that $0\leq j\leq k$%
\begin{equation}
\mathbb{E}(Y_{0}Y_{2j})\geq\mathbb{E}(Y_{0}Y_{2k})\geq0, \label{l3}%
\end{equation}
and for integers $k\geq2$%
\begin{equation}
\mathbb{E}(Y_{0}Y_{k})\leq\mathbb{E}(Y_{0}Y_{2}). \label{l4}%
\end{equation}
For any positive integer $\ell$ and any $j\geq2\ell$%
\begin{equation}
\sum_{k=2\ell}^{j}\mathbb{E}(Y_{0}Y_{k})\geq0. \label{l5}%
\end{equation}
For any positive integer $\ell$ and any $n\geq2\ell$
\begin{equation}
\max_{2\ell\leq j\leq n}\sum_{k=2\ell}^{j}\mathbb{E}(Y_{0}Y_{k})\leq
\sum_{k=2\ell}^{n}\mathbb{E}(Y_{0}Y_{k})+\mathbb{E}(Y_{0}Y_{2\ell}).
\label{l6}%
\end{equation}

\end{lemma}

\textbf{Proof. }To prove relation (\ref{l1}) just note that $\mathbb{E}%
(Y_{0}Y_{k+j})=\mathbb{E}(Y_{0}\mathbb{E}_{0}Y_{k+j})$ and apply relation
(\ref{spectral}). Then, note that%
\[
\mathbb{E}(Y_{0}(Y_{2k}+Y_{2k+1}))=\int_{-1}^{1}(s^{2k}+s^{2k+1}%
)\nu(\mathrm{d}s),
\]
and $s^{2k}(1+s)\geq0$ for all $-1\leq s\leq1;$ so relation (\ref{l2}) holds
by (\ref{spectral}).  Relation (\ref{l3}) is clearly true since for $0\leq
j\leq k,$ we have $s^{2j}\geq$ $s^{2k}\geq0$ for all $s$. Finally in order to
show relation (\ref{l4}) just note that for all $-1\leq s\leq1$ and any
integer $k\geq2$ we have $s^{k}\leq s^{2},$ which we combine with
(\ref{spectral}).

Relations (\ref{l5}) and (\ref{l6}) are obtained via a blocking argument. Now,
if $\ j$ is odd, say $j=2m+1$, we can write%
\[
\sum_{k=2\ell}^{j=2m+1}\mathbb{E}(Y_{0}Y_{k})=\sum_{k=\ell}^{m}(\mathbb{E}%
(Y_{0}Y_{2k})+\mathbb{E}(Y_{0}Y_{2k+1}))\geq0,
\]
since, by relation (\ref{l2}) in the right hand side we have a sum of positive terms.

On the other hand, if $j$ is even, say $j=2(m+1)$ \
\[
\sum_{k=2\ell}^{j=2m+2}\mathbb{E}(Y_{0}Y_{k})=\mathbb{E}(Y_{0}Y_{2(m+1)}%
)+\sum_{k=2\ell}^{2m+1}\mathbb{E}(Y_{0}Y_{k}).
\]
By relation (\ref{l3})$,$ $\mathbb{E}(Y_{0}Y_{2(m+1)})\geq0$. Therefore
(\ref{l5}) is true for all $j\geq2\ell.$

Now, if $n\geq2\ell,$ by combining the latter considerations with (\ref{l3})
and (\ref{l4})%
\begin{gather*}
\max_{2\ell\leq j\leq n}\sum_{k=2\ell}^{j}\mathbb{E}(Y_{0}Y_{k})\leq
\max_{2\ell\leq2j\leq n}\sum_{k=2\ell}^{2j}\mathbb{E}(Y_{0}Y_{k})\vee
\max_{2\ell\leq2m+1\leq n}\sum_{k=2\ell}^{2m+1}\mathbb{E}(Y_{0}Y_{k})\\
\leq\max_{2\ell\leq2m+1\leq n}\sum_{k=2\ell}^{2m+1}\mathbb{E}(Y_{0}Y_{k}%
)+\max_{2\ell\leq2j\leq n}\mathbb{E}(Y_{0}Y_{2j})\leq\mathbb{E}(Y_{0}Y_{2\ell
})+\sum_{k=2\ell}^{m_{n}}\mathbb{E}(Y_{0}Y_{k}).
\end{gather*}
where $m_{n}$ is the largest odd integer smaller than $n.$ If $n$ is odd
$n=m_{n}.$ If $n$ is even, $n>2\ell,$ then $m_{n}=n-1.$ By taking into account
relation (\ref{l3}) we can add in this case a positive term, $\mathbb{E}%
(Y_{0}Y_{n}),$ and obtain overall relation (\ref{l6}). $\ \square$

\bigskip

Next, we give a CLT for a triangular array of row-wise stationary reversible
Markov chains. The conditions for the CLT are imposed to the covariances of
both the variables and their squares$.$

\begin{theorem}
\label{CLT-tri-rev} Let $(X_{i})_{i\in\mathbb{Z}}$ be a stationary reversible
Markov chain. For real functions $f_{n},$ define
\begin{equation}
X_{n,k}=f_{n}(X_{k}). \label{defX}%
\end{equation}
Assume that
\begin{equation}
\mathbb{E}X_{n,k}^{4}<\infty;~\mathbb{E}X_{n,k}=0\text{ and }\mathbb{E}%
(X_{n,0}^{2})\rightarrow\sigma^{2}, \label{firstcond}%
\end{equation}%
\begin{equation}
\mathrm{cov}(X_{n,0},X_{n,2})+\sum_{k=2}^{n}\mathrm{cov}(X_{n,0}%
,X_{n,k})\rightarrow0, \label{neglcov}%
\end{equation}
and
\begin{equation}
\frac{1}{n}(\mathrm{var}(X_{n,0}^{2})+\sum_{u=0}^{n}\mathrm{cov}(X_{n,0}%
^{2},X_{n,u}^{2}))\rightarrow0. \label{secondcond}%
\end{equation}
Then
\[
\frac{1}{\sqrt{n}}\sum_{k=1}^{n}X_{n,k}\Rightarrow N(0,\sigma^{2}).
\]

\end{theorem}

\textbf{Proof}. We start from a standard martingale decomposition by using
projections:
\begin{align*}
S_{n}  &  =\sum_{k=1}^{n}(X_{n,k}-\mathbb{E}_{k-1}X_{n,k})+\sum_{k=1}%
^{n}\mathbb{E}_{k-1}X_{n,k}\\
&  =\sum_{k=1}^{n}D_{n,k}+\sum_{k=1}^{n}\mathbb{E}_{k-1}X_{n,k}.
\end{align*}
Note that $D_{n,k}=X_{n,k}-\mathbb{E}_{k-1}X_{n,k}$ are martingale differences
adapted to $(\mathcal{F}_{k}\mathcal{)}_{k\geq1}.$

We show first that the second term divided by $\sqrt{n}$ is negligible for the
convergence in distribution. To show this we estimate $\mathrm{var}(\sum
_{k=1}^{n}\mathbb{E}_{k-1}X_{n,k}).$ By using the properties of conditional
expectation, stationarity and relation (\ref{l1}) in Lemma \ref{cov}, for
$k\leq j$ we obtain%

\begin{align*}
\mathbb{E}(\mathbb{E}_{k-1}X_{n,k}\mathbb{E}_{j-1}X_{n,j})  &  =\mathbb{E}%
(X_{n,j}\mathbb{E}_{k-1}X_{n,k})=\mathbb{E}(X_{n,j-k+1}\mathbb{E}_{0}%
X_{n,1})\\
&  =\mathbb{E}(\mathbb{E}_{0}X_{n,1}\mathbb{E}_{0}X_{n,j-k+1})=\mathbb{E}%
(X_{n,0}X_{n,j-k+2})
\end{align*}
and so,
\begin{gather*}
\mathbb{E}(\sum_{k=1}^{n}\mathbb{E}_{k-1}X_{n,k})^{2}=n\cdot cov(X_{n,0}%
,X_{n,2})+2\sum_{j=2}^{n}\sum_{k=1}^{j-1}\mathbb{E}(X_{n,0}X_{n,j-k+2})\\
\leq2n\max_{2\leq j\leq n}\sum_{k=2}^{j}\mathbb{E}(X_{n,0}X_{n,k}).
\end{gather*}
By applying now relation (\ref{l6}) of Lemma \ref{cov} $\ $%
\[
\frac{1}{n}\mathbb{E}(\sum_{k=1}^{n}\mathbb{E}_{k-1}X_{n,k})^{2}%
\leq2\mathbb{E}(X_{n,0}X_{n,2})+2\sum_{k=2}^{n}\mathbb{E}(X_{n,0}X_{n,k}),
\]
which converges to $0$ by (\ref{neglcov}).

Now we analyze the martingale differences via Theorem \ref{Thmart} given in
Appendix. To show that $\max_{1\leq k\leq n}|D_{n,k}|/\sqrt{n}$ is uniformly
integrable we show that $\mathbb{E}(\max_{1\leq k\leq n}D_{n,k}^{2})\leq Cn,$
for all $n$ and some constant $C>0$. Indeed, since $\mathbb{E}(D_{n,0}%
^{2})\leq\mathbb{E}(X_{n,0}^{2})$, by (\ref{firstcond}) we note that there is
a positive constant $C$ such that $\mathbb{E}(D_{n,0}^{2})\leq C$ and
therefore, by stationarity%
\[
\mathbb{E}(\max_{1\leq k\leq n}D_{n,k}^{2})\leq\sum_{k=1}^{n}\mathbb{E}%
(D_{n,k}^{2})\leq Cn.
\]
It remains to verify
\[
\frac{1}{n}\sum_{k=1}^{[nt]}D_{n,k}^{2}\xrightarrow{P}t\sigma^{2}.
\]
We note that%
\[
D_{n,k}^{2}=X_{n,k}^{2}+(\mathbb{E}_{k-1}X_{n,k})^{2}-2X_{n,k}(\mathbb{E}%
_{k-1}X_{n,k})=X_{n,k}^{2}+I_{n,k}.
\]
Furthermore, by the Cauchy-Schwartz inequality, (\ref{firstcond}), Lemma
\ref{cov} and stationarity
\begin{gather*}
\frac{1}{n}\mathbb{E}\sum_{k=1}^{n}|I_{n,k}|=\frac{1}{n}\mathbb{E}\sum
_{k=1}^{n}|(\mathbb{E}_{k-1}X_{n,k})^{2}-2X_{n,k}(\mathbb{E}_{k-1}X_{n,k})|\\
\leq\mathbb{E}((\mathbb{E}_{-1}X_{n,0})^{2})+\frac{2}{n}\sum_{k=1}%
^{n}||X_{n,k}||_{2}||\mathbb{E}_{k-1}X_{n,k}||_{2}\\
\leq\mathrm{cov}(X_{n,0},X_{n,2})+C\sqrt{\mathrm{cov}(X_{n,0},X_{n,2})}.
\end{gather*}
We see that the last quantity converges to $0$ by condition (\ref{neglcov})
combined with relation (\ref{l5}) in Lemma \ref{cov}. We also note that by
stationarity and (\ref{firstcond}),
\[
\frac{1}{n}\mathbb{E}(\sum_{k=1}^{[nt]}X_{n,k}^{2})=\frac{[nt]}{n}%
\mathbb{E}(X_{n,0}^{2})\rightarrow\sigma^{2}t.
\]
So, it remains\ to show that%
\[
\frac{1}{n}\sum_{k=1}^{[nt]}(X_{n,k}^{2}-\mathbb{E}(X_{n,0}^{2}%
))\xrightarrow{P}0\text{,}%
\]
which will be implied by
\[
\mathrm{var}(\frac{1}{n}\sum_{k=1}^{[nt]}X_{n,k}^{2})\rightarrow0.
\]
We estimate now this variance. Note that by relation (\ref{l6}) in Lemma
\ref{cov},
\begin{gather*}
\frac{1}{n^{2}}\sum_{k=1}^{[nt]}\mathrm{var}(X_{n,k}^{2})+\frac{2}{n^{2}}%
\sum_{k=2}^{[nt]}\sum_{u=1}^{k-1}\mathrm{cov}(X_{n,k}^{2},X_{n,u}^{2}%
)\leq\frac{2[nt]}{n^{2}}(\max_{0\leq k\leq n}\sum_{u=0}^{k}\mathrm{cov}%
(X_{n,0}^{2},X_{n,u}^{2}))\\
\leq\frac{2t}{n}(\mathrm{var}(X_{n,0}^{2})+\sum_{u=0}^{n}\mathrm{cov}%
(X_{n,0}^{2},X_{n,u}^{2})).
\end{gather*}
The result follows by condition (\ref{secondcond}). $\ \square$

\section{Proof of Theorem\textbf{ \ref{gclt}}}

By the consistency of $\hat{f}_{n}(x_{j})$ due to the continuity of $f(x)$ at
$x_{j}$ and the assumptions we made on the bandwidth and kernel we only need
to show that
\[
\sqrt{nb_{n}}\left(  \frac{\hat{f}_{n}(x_{j})-\mathbb{E}\hat{f}_{n}(x_{j}%
)}{(f(x_{j})\int K^{2}(u)du)^{1/2}},1\leq j\leq m\right)  \Rightarrow
N(0,I_{m}).
\]
By the Cramer-Wold device, it suffices to prove that%
\[
\sqrt{nb_{n}}\sum_{j=1}^{m}\lambda_{j}\frac{\hat{f}_{n}(x_{j})-\mathbb{E}%
\hat{f}_{n}(x_{j})}{(f(x_{j})\int K^{2}(u)du)^{1/2}}\Rightarrow\sum_{j=1}%
^{m}\lambda_{j}Z_{j}%
\]
for arbitrary fixed $\lambda_{1},\cdots,\lambda_{m}$. Here $Z_{j},1\leq j\leq
m$ are i.i.d. standard normal random variables.

Let $S_{n}=\sum_{i=1}^{n}Y_{n,i}$ where%

\[
Y_{n,i}=\frac{1}{\sqrt{b_{n}}}\sum_{j=1}^{m}\frac{\lambda_{j}}{(f(x_{j})\int
K^{2}(u)du)^{1/2}}(K(\frac{x_{j}-X_{i}}{b_{n}})-\mathbb{E}K(\frac{x_{j}-X_{i}%
}{b_{n}})).
\]
We shall verify the conditions of Theorem \ref{CLT-tri-rev}. We verify first
(\ref{firstcond}).
\begin{gather}
\mathrm{var}(Y_{n,i})=\frac{1}{b_{n}}\sum_{j=1}^{m}\frac{\lambda_{j}^{2}%
}{f(x_{j})\int K^{2}(u)du}\mathrm{var}(K(\frac{x_{j}-X_{0}}{b_{n}%
}))\nonumber\\
+\frac{2}{b_{n}}\sum_{p=1}^{m-1}\sum_{j=p+1}^{m}\frac{\lambda_{j}\lambda_{p}%
}{(f(x_{j})f(x_{p}))^{1/2}\int K^{2}(u)du}\mathrm{cov}(K(\frac{x_{j}-X_{0}%
}{b_{n}}),K(\frac{x_{p}-X_{0}}{b_{n}}))=I_{n}+II_{n}. \label{cross}%
\end{gather}
Now, by Bochner's lemma (see Parzen, 1962, or Bosq, 1998) and the fact
$b_{n}\rightarrow0,$%
\[
\lim_{n\rightarrow\infty}\mathrm{var}(\frac{1}{\sqrt{b_{n}}}K(\frac
{x_{j}-X_{0}}{b_{n}}))=f(x_{j})\int K^{2}(u)du.
\]
Therefore
\[
I_{n}\rightarrow\sum_{j=1}^{m}\lambda_{j}^{2}.
\]
On the other hand by simple calculus computations involving the symmetry of
$K,$ for $j\neq p$ we have%
\begin{align*}
&  \frac{1}{b_{n}}\mathrm{cov}(K(\frac{x_{j}-X_{0}}{b_{n}}),K(\frac
{x_{p}-X_{0}}{b_{n}}))\\
&  =\frac{1}{b_{n}}\int K(\frac{x_{j}-u}{b_{n}})K(\frac{x_{p}-u}{b_{n}%
})f(u)du-\frac{1}{b_{n}}\int K(\frac{x_{j}-u}{b_{n}})f(u)du\int K(\frac
{x_{p}-u}{b_{n}})f(u)du\\
&  =\int K(v)K(v+\frac{x_{p}-x_{j}}{b_{n}})f(x_{j}-b_{n}v)dv-b_{n}\int
K(v)f(x_{j}-b_{n}v)dv\int K(v)f(x_{p}-b_{n}v)dv.
\end{align*}
Clearly the second term is convergent to $0$ by Bochner's lemma and the fact
that $b_{n}\rightarrow0.$ For the first term we cannot apply directly the
Bochner lemma, but by using the same arguments as in its proof presented in
(Parzen, 1962) along with the Lebesgue dominated convergence theorem, under
our conditions we deduce that this term is also negligible. Hence
$Var(Y_{n,i})\rightarrow\sum_{j=1}^{m}\lambda_{j}^{2}$.

To verify condition (\ref{secondcond}) we introduce the function
\[
\tilde{g}(u)=\left(  \sum_{j=1}^{m}\frac{\lambda_{j}}{(f(x_{j})\int
K^{2}(u)du)^{1/2}}(K(\frac{x_{j}-u}{b_{n}})-\mathbb{E}K(\frac{x_{j}-X_{0}%
}{b_{n}}))\right)  ^{2}.
\]
Since by our conditions on $K,$ the function $\tilde{g}(u)$ has bounded
derivative, by Newman extension of Hoeffding lemma (see relation (22) in
Newman, 1980),
\[
\mathrm{cov}(Y_{n,0}^{2},Y_{n,k}^{2})=\frac{1}{b_{n}^{2}}\int\int\tilde
{g}^{\prime}(u)\tilde{g}^{\prime}(v)H_{k}(u,v)dudv.
\]
Therefore, with $C=(\sum_{j=1}^{m}\lambda_{j}/(f(x_{j})\int K^{2}%
(u)du)^{1/2})^{4}$ we obtain
\begin{equation}
\frac{1}{n}\sum_{k=1}^{n}|cov(Y_{n,0}^{2},Y_{n,k}^{2})|\ll\frac{16C}%
{nb_{n}^{4}}||KK^{\prime}||_{\infty}^{2}\sum_{k=1}^{n}\eta_{k},\nonumber
\end{equation}
which converges to $0$ by taking into account our conditions on $b_{n}$ and
$\eta_{k}$.

We have also to treat $\mathrm{var}(Y_{n,0}^{2})/n.$ We shall apply first
H\"{o}lder inequality to obtain%

\[
\mathrm{var}(Y_{n,0}^{2})\leq\mathbb{E}(Y_{n,0}^{4})\leq\frac{8m^{3}}%
{b_{n}^{2}}\sum_{j=1}^{m}\frac{\lambda_{j}^{4}}{(f(x_{j})\int K^{2}(u)du)^{2}%
}(K^{4}(\frac{x_{j}-X_{0}}{b_{n}})+(\mathbb{E}K(\frac{x_{j}-X_{0}}{b_{n}%
}))^{4}).
\]
Note that for any $p\geq1$%
\[
\mathbb{E}K^{p}(\frac{x-X_{0}}{b_{n}})=b_{n}\int K^{p}(u)f(x+b_{n}u)du.
\]
So by the Bochner's lemma and the fact that the kernel is bounded%
\[
\frac{1}{n}\mathrm{var}(Y_{n,0}^{2})\leq\max_{1\leq j\leq m}\frac{C_{m}%
}{nb_{n}}\int K^{4}(u)f(x_{j}+b_{n}u)du+\frac{4}{n}b_{n}^{2}(\int
K(u)f(x_{j}+b_{n}u)du)^{4}\rightarrow0
\]
provided $nb_{n}\rightarrow\infty$. Therefore (\ref{secondcond}) is satisfied.

We turn now to verify condition (\ref{neglcov}). Since%
\[
\mathrm{cov}(Y_{n,0},Y_{n,k})=\frac{1}{b_{n}}\sum_{j,p=1}^{m}\frac{\lambda
_{j}\lambda_{p}}{(f(x_{j})f(x_{p}))^{1/2}\int K^{2}(u)du}\mathrm{cov}%
(K(\frac{x_{j}-X_{0}}{b_{n}}),K(\frac{x_{p}-X_{k}}{b_{n}})),
\]
it is enough to show that for any $j$ and $p$ fixed%
\begin{equation}
\frac{1}{b_{n}}\mathrm{cov}(K(\frac{x_{j}-X_{0}}{b_{n}}),K(\frac{x_{p}-X_{2}%
}{b_{n}}))\rightarrow0 \label{negl1}%
\end{equation}
and
\begin{equation}
\sum_{k=2}^{n}\frac{1}{b_{n}}\mathrm{cov}(K(\frac{x_{j}-X_{0}}{b_{n}}%
),K(\frac{x_{p}-X_{k}}{b_{n}}))\rightarrow0. \label{negl2}%
\end{equation}
We shall estimate $\mathrm{cov}(X_{n,0},X_{n,k})$ in two different ways and
take the minimum of these estimates.

By Lemma \ref{cov}, for $k\geq2$%
\begin{gather*}
A=\mathrm{cov}(K(\frac{x_{j}-X_{0}}{b_{n}}),K(\frac{x_{p}-X_{k}}{b_{n}}%
))\leq\mathrm{cov}(K(\frac{x_{j}-X_{0}}{b_{n}}),K(\frac{x_{p}-X_{2}}{b_{n}%
}))\\
=\frac{1}{b_{n}}\int\int K(\frac{x_{j}-u}{b_{n}})K(\frac{x_{p}-v}{b_{n}%
})(f_{2}(u,v)-f(u)f(v))dudv.
\end{gather*}
By changing the variable
\begin{gather}
A=b_{n}\int\int K(u)K(v)(f_{2}(x_{j}+ub_{n},x_{p}+vb_{n})-f(x_{j}%
+ub_{n})f(x_{p}+vb_{n}))dudv\label{formula}\\
\leq b_{n}\int\int K(u)K(v)(f_{2}(x_{j}+ub_{n},x_{p}+vb_{n})dudv.\nonumber
\end{gather}
To analyze this term we divide the integral in (\ref{formula}) on $4$ sets

$(|u-x_{j}|\leq M)\times(|v-x_{p}|\leq M)$

$(|u-x_{j}|>M)\times(|v-x_{p}|\leq M)$

$(|u-x_{j}|\leq M)\times(|v-x_{p}|>M)$

$(|u-x_{j}|\geq\ M)\times(|v-x_{p}|>M)$.\ 

On the first set, $(|u-x_{j}|\leq M)\times(|v-x_{p}|\leq M),$ we change the
variable and obtain%
\begin{gather*}
b_{n}\int_{-\frac{M}{b_{n}}}^{\frac{M}{b_{n}}}\int_{-\frac{M}{b_{n}}}%
^{\frac{M}{b_{n}}}K(u)K(v)(f_{2}(x_{j}-b_{n}u,x_{p}-b_{n}v)dudv\leq\\
\leq b_{n}^{\ }\sup_{|a|<M}|f_{2}(x_{j}+a,x_{p}+a)|\int\int K(u)K(v)dudv.
\end{gather*}
By our assumptions this term is smaller than $b_{n}^{\ }C_{M}.$ On the set
$(|u-x_{j}|>M)\times(|v-x_{p}|\leq M)$ we have
\[
\frac{1}{b_{n}}\mathbb{E}(K(\frac{x_{j}-X_{0}}{b_{n}})K(\frac{x_{p}-X_{2}%
}{b_{n}})I(|X_{0}-x_{j}|>M)I(|X_{2}-x_{p}|\leq M)\leq||K||_{\infty}\frac
{1}{b_{n}}K(\frac{M}{b_{n}}).
\]
A similar estimate is obtained on the set $(|u-x_{j}|\leq M)\times
(|v-x_{p}|>M).$ On $(|u-x_{j}|>M)\times(|v-x_{p}|>M)$ we estimate in the
following way%
\[
\frac{1}{b_{n}}\mathbb{E}(K(\frac{x_{j}-X_{0}}{b_{n}})K(\frac{x_{p}-X_{2}%
}{b_{n}})I(|X_{0}-x_{j}|>M)I(|X_{2}-x_{p}|>M)\leq\frac{1}{b_{n}}K^{2}(\frac
{M}{b_{n}}).
\]
So, for $n$ sufficiently large
\begin{equation}
\frac{1}{b_{n}}\mathrm{cov}(K(\frac{x_{j}-X_{0}}{b_{n}}),K(\frac{x_{p}-X_{2}%
}{b_{n}}))\leq b_{n}C(x_{j},x_{p}) \label{rel1}%
\end{equation}
where%
\[
C(x_{j},x_{p})=C_{M}+2||K||_{\infty}\frac{1}{b_{n}^{2}}K(\frac{M}{b_{n}%
})+\frac{1}{b_{n}^{2}}K^{2}(\frac{M}{b_{n}}).
\]
By our conditions on $K$ for $n$ sufficiently large $C(x_{j},x_{p})$ is bounded.

So clearly by (\ref{rel1}) condition (\ref{negl1}) is satisfied.

On the other hand,
\begin{gather}
\frac{1}{b_{n}}\mathrm{cov}(K(\frac{x_{j}-X_{0}}{b_{n}}),K(\frac{x_{p}-X_{k}%
}{b_{n}}))\label{rel2}\\
=\frac{1}{b_{n}^{3}}\int%
{\displaystyle\int}
K^{\prime}(\frac{x_{j}-u}{b_{n}})K^{\prime}(\frac{x_{p}-v}{b_{n}}%
)H_{k}(u,v)dudv\leq||K^{\prime}||_{\infty}^{2}\frac{1}{b_{n}^{3}}\eta
_{k}.\nonumber
\end{gather}
So by combining the estimates in (\ref{rel1}) and (\ref{rel2}) we have proven
that
\begin{align}
&  \frac{1}{b_{n}}\mathrm{cov}(K(\frac{x_{j}-X_{0}}{b_{n}}),K(\frac
{x_{p}-X_{k}}{b_{n}})\label{rel3}\\
&  \leq\min(b_{n}C(x_{j},x_{p}),c\frac{1}{b_{n}^{3}}cov(X_{0},X_{k}))\ll
\min(b_{n},\frac{1}{b_{n}^{3}}\eta_{k}).\nonumber
\end{align}
To continue we use the estimate from (\ref{rel3}) to bound the sum in the
right hand side of (\ref{negl2}). We shall divide the sum in two, up to
$m_{n}$ and after $m_{n}.$ This sequence of positive integers $m_{n}$ will be
selected later. On the first part of the sum we use the bound of order $b_{n}$
and on the second part of the sum we use the bound $\eta_{k}/b_{n}^{3}$. So,
by the properties of slowly varying function $l,$
\begin{align*}
\frac{1}{b_{n}}\sum_{k=2}^{n}\mathrm{cov}(K(\frac{x_{j}-X_{0}}{b_{n}}%
),K(\frac{x_{p}-X_{k}}{b_{n}}))  &  \ll\sum_{k=2}^{m_{n}}b_{n}+\sum
_{k=m_{n}+1}^{n}\frac{1}{b_{n}^{3}}\eta_{k}\\
&  \ll m_{n}b_{n}+\frac{1}{b_{n}^{3}m_{n}^{3}l\left(  m_{n}\right)  }.
\end{align*}
To optimize the sum we take
\[
m_{n}=[\max(\frac{1}{\sqrt{b_{n}}},\text{ }\frac{1}{b_{n}l^{1/6}(1/b_{n}%
^{1/2})})]+1.
\]
Clearly $m_{n}b_{n}\rightarrow0\ ($since $b_{n}/\sqrt{b_{n}}\rightarrow0$ and
$1/l^{1/6}(1/b_{n}^{1/2})\rightarrow0).$

Since $m_{n}>b_{n}^{-1/2}$ and $l$ is increasing we have $l(m_{n}%
)>l(b_{n}^{-1/2})$ and so, since $m_{n}>(b_{n}l^{1/6}(1/b_{n}^{1/2}))^{-1},$
\[
\frac{1}{b_{n}^{3}m_{n}^{3}l_{m_{n}}}\leq\frac{1}{b_{n}^{3}m_{n}^{3}%
l(b_{n}^{-1/2})}\leq\frac{l^{1/2}(1/b_{n}^{1/2})}{l(b_{n}^{-1/2})}\leq\frac
{1}{\sqrt{l(1/b_{n}^{1/2})}}\rightarrow0\text{ as }n\rightarrow\infty.
\]
By taking now into account Lemma \ref{cov}, it follows that the sum in
(\ref{negl2}) is positive and then (\ref{negl2})~follows. Now by (\ref{negl1})
and (\ref{negl2}) we conclude that condition (\ref{neglcov}) is satisfied and
the result follows. $\ \square$

\bigskip

\textbf{Proof of Remark \ref{rmclt}. }To prove this remark we have to replace
in the proof of Theorem \ref{gclt} relation (\ref{rel2}) by relation (3.12) in
Bosq et al. (1999), namely%
\[
\frac{1}{b_{n}}|\mathrm{cov}(K(\frac{x_{j}-X_{0}}{b_{n}}),K(\frac{x_{p}-X_{k}%
}{b_{n}})|\leq\frac{4}{b_{n}}||K||_{\infty}^{2}\alpha_{k}.
\]
Then, we replace relation (\ref{rel3}) by
\[
\frac{1}{b_{n}}\mathrm{cov}(K(\frac{x_{j}-X_{0}}{b_{n}}),K(\frac{x_{p}-X_{k}%
}{b_{n}})\ll\min(b_{n},\frac{1}{b_{n}}\alpha_{k}),
\]
and follow the proof from the page 88 in Bosq et al. (1999), to obtain the
result of this remark.

\section{Appendix}

\textbf{Martingale limit theorem} (G\"{a}nssler and H\"{a}usler, 1986, pages 315-317).

\begin{theorem}
\label{Thmart}Assume $(D_{n,k})_{1\leq k\leq n}$ is a triangular array of
martingales adapted to an increasing in $k$ filtration $\mathcal{F}_{n,k}.$
Assume $\sum_{k=1}^{[nt]}D_{n,k}^{2}\xrightarrow{P}t\sigma^{2}$ and
\begin{equation}
\max_{1\leq k\leq n}|D_{n,k}|\text{ is uniformly integrable} \label{cond2}%
\end{equation}
(as before, by $[x]$ we denote as usual the integer part of $x$). Then
$S_{[nt]}\Rightarrow\sigma W(t)$ where $W(t)$ is a standard Brownian measure.
In particular $S_{n}\Rightarrow N(0,\sigma^{2}).$
\end{theorem}

\bigskip

\textbf{Acknowledgement.} The authors would like to thank the referee for
numerous suggestions which improved the presentation of this paper.

\end{document}